\documentclass[3p]{elsarticle}

\usepackage{amsmath}
\usepackage{amsthm}
\usepackage{amssymb}
\usepackage{bm}
\usepackage{accents}
\usepackage{color}

\newdefinition{rmk}{Remark}
\newproof{pf}{Proof}

\numberwithin{equation}{section}

\graphicspath{{figs/}}

\journal{arXiv}

\begin{document}

\begin{frontmatter}

\title{Numerical solution of 2D boundary value problems on merged Voronoi--Delaunay meshes\tnoteref{label1}}
\tnotetext[label1]{The work was supported by the Russian Science Foundation (grant No. 24-11-00058).}

\author{M. M. Chernyshov\fnref{lab1}}
\ead{chernyshovmm@my.msu.ru}

\author{P.N. Vabishchevich\corref{cor1}\fnref{lab1,lab2}}
\ead{vab@cs.msu.ru}

\cortext[cor1]{Correspondibg author.}

\address[lab1]{Lomonosov Moscow State University, 1, building 52, Leninskie Gory,  119991 Moscow, Russia}

\address[lab2]{North-Eastern Federal University, 58, Belinskogo st, Yakutsk, 677000, Russia}

\begin{abstract}

Computational technologies for the approximate solution of multidimensional boundary value problems often rely on irregular computational meshes and finite-volume approximations.
In this framework, the discrete problem represents the corresponding conservation law for control volumes associated with the nodes of the mesh.
This approach is most naturally and consistently implemented using Delaunay triangulations together with Voronoi diagrams as control volumes.
In this paper, we employ meshes with nodes located both at the vertices of Delaunay triangulations and at the generators of Voronoi partitions.
The cells of the merged Voronoi--Delaunay mesh are orthodiagonal quadrilaterals.
On such meshes, scalar and vector functions, as well as invariant gradient and divergence operators of vector calculus, can be conveniently approximated.
We illustrate the capabilities of this approach by solving a steady--state diffusion--reaction problem in an anisotropic medium.
\end{abstract}

\begin{keyword}
Elliptic BVP \sep finite volume method \sep Delaunay triangulation \sep Voronoi tessellation \sep discrete vector operators
\MSC  35J25 \sep 65N08 \sep 65D18 \sep 65N06
\end{keyword}
\end{frontmatter}

\section{Introduction}

Applied mathematical models provide quantitative tools for analyzing and predicting physical, chemical, and biological processes \cite{LionsBook}.
Their theoretical description is based on systems of coupled partial differential equations, which include elliptic, parabolic, and hyperbolic types \cite{evans2010partial}.
The corresponding initial--boundary value problems are formulated in multidimensional domains with complex geometry, and numerical methods are the primary tool for their solution \cite{KnabnerAngermann2003,QuarteroniValli}.

Modern computational practice makes extensive use of unstructured meshes for irregular geometries \cite{thompson1998handbook}.
For two-dimensional problems, triangular meshes are most common, while in three dimensions tetrahedral meshes are used \cite{lo2014finite}.
Spatial approximations are most often constructed using the finite element method (FEM) \cite{pdeinc2009solutions,brenner2008mathematical}.
The widely known FEniCS computational platform \cite{AlnaesEtal2015}, for example, is based on FEM technology.

Finite difference methods for boundary value problems are well developed on regular meshes \cite{Samarskii1989,Strikwerda2004}.
On irregular meshes, approximate solutions are constructed using the balance (integro--interpolation) method \cite{TikhonovSamarskii1961}.
In this case, the discrete problem is regarded as a conservation law for individual cells of the computational domain (control volumes).
This approach is now widely referred to as the finite volume method (FVM) \cite{FVM,leveque2002finite}.
A finite-element version of FVM is presented in \cite{Generalized_difference_methods}.

A key element of finite-volume technology is the use of irregular meshes, where the approximate solution is defined at the nodes.
For each mesh node, a control volume is constructed, and a balance relation is written, which formulates the conservation law at the discrete level.
This requires the introduction of flux quantities that provide interface relations at the boundaries of control volumes.
Such discrete problems can be related to the simplest approximations using mixed finite elements \cite{boffi2013mixed}.

Finite-volume approximations are most straightforwardly built on optimal meshes: Delaunay triangulations in two dimensions and their tetrahedral analogues in three dimensions \cite{george1998delaunay,cheng2013delaunay}.
In this case, the control volumes are Voronoi polygons.
It is crucial that the mesh operators of vector calculus, invariant under a change of coordinates, inherit the main properties of the corresponding differential operators by satisfying the basic integral relations \cite{shashkov2018conservative}.
Such consistent approximations form the basis of mimetic discretization (MD) technology \cite{da2014mimetic,vabishchevich2005finite}.

In standard mesh generation, the nodes of the computational mesh are placed at the vertices of the cells.
Another possibility is to define the solution at the cell centers (e.g., at their centroids).
In this case, the cell itself serves as a natural control volume for finite-volume approximation.
In a Delaunay triangulation, Voronoi polyhedra uniquely define control volumes if the solution is assigned to mesh nodes.
Conversely, if the solution is defined at the vertices of the Voronoi diagram, then the control volumes coincide with the cells of the Delaunay mesh.
This duality motivates the use of both Delaunay and Voronoi meshes in finite-volume approximations, thereby extending the class of solvable problems.

The idea of a two-mesh technology employing both Delaunay and Voronoi meshes was proposed earlier in \cite{frjazinov1975,samarskiiUMN}.
There are still relatively few works in this direction.
The recent study \cite{vabishchevich2024gvm} constructs operator-difference approximations on a merged rectangular mesh for scalar and vector boundary value problems.
A general treatment of vector operators on merged two-dimensional Voronoi--Delaunay meshes is presented in \cite{MVD-2d}.
Here the computational domain is partitioned into orthodiagonal quadrilaterals, and consistent approximations for gradient, divergence, and curl operators are developed.
The applicability of this approach to mesh formulations of scalar and vector boundary value problems is briefly discussed.

In this work, we construct finite-volume approximations for a two-dimensional diffusion--reaction problem in an anisotropic medium.
In this setting, the use of a single mesh (Delaunay or Voronoi) is problematic, while the merged Voronoi--Delaunay mesh offers an effective alternative.
We emphasize the main features of the approach and present numerical results for a model problem, which demonstrate the efficiency of FVM on such meshes.

The structure of the paper is as follows.
Section~2 formulates the two-dimensional boundary value problem for the diffusion--reaction equation with anisotropic coefficients.
Section~3 describes Delaunay and Voronoi meshes, scalar and vector mesh functions, and introduces the merged quadrilateral Voronoi--Delaunay mesh.
Section~4 develops approximations of the gradient and divergence operators.
Section~5 presents the mesh boundary value problem and establishes an a priori stability estimate.
Section~6 provides numerical results for test problems with different diffusion tensors.
Section~7 concludes with a summary of the main findings.

\section{Problem formulation}

We consider the linear steady--state diffusion--reaction problem in a bounded domain $\Omega \subset \mathbb{R}^2$.
The concentration of the substance $u(\bm x)$, $\bm x = (x^{(1)}, x^{(2)})$, satisfies the equation
\begin{equation}\label{2.1}
	- \operatorname{div} \!\left( K(\bm x)\, \operatorname{grad} u(\bm x) \right) + r(\bm x)\, u(\bm x) = f(\bm x), 
  \quad \bm x \in \Omega .
\end{equation}
Equation \eqref{2.1} is supplemented with homogeneous Dirichlet boundary conditions:
\begin{equation}\label{2.2}
	u(\bm x) = 0, \quad \bm x \in \partial\Omega .
\end{equation}

The medium is assumed to be anisotropic.  
The diffusion coefficient $K(\bm x) = \{k_{\alpha \beta}(\bm x)\}$ is a symmetric, positive-definite tensor of rank two:
\[
  k_{\alpha \beta}(\bm x) = k_{\beta \alpha}(\bm x), \quad 
  \alpha, \beta = 1,2, \quad
  \sum_{\alpha,\beta=1}^{2} k_{\alpha \beta}(\bm x)\, \zeta_{\alpha} \zeta_{\beta} \ge \kappa \sum_{\alpha=1}^{2} \zeta_{\alpha}^2,
  \quad \kappa > 0, \quad \bm x \in \Omega ,
\]
for any vector $\bm \zeta = (\zeta_1,\zeta_2)$.  
The reaction coefficient $r(\bm x)$ is assumed nonnegative: $r(\bm x) \ge 0$.

Let $\mathcal{H} = L_2(\Omega)$ with inner product and norm
\[
  (v,w) = \int_{\Omega} v(\bm x)\, w(\bm x)\, d\bm x, 
  \quad \|v\| = (v,v)^{1/2}.
\]
For solutions of problem \eqref{2.1}, \eqref{2.2}, one can derive a number of a priori estimates. In particular, the fundamental estimate
\begin{equation}\label{2.3}
  \|\operatorname{grad} u\| \le \frac{c}{\kappa}\, \|f\|
\end{equation}
holds, where $c$ is the constant from Friedrichs' inequality.  
This relation will serve as a reference when analyzing discrete problems.

To construct finite-volume approximations of problem \eqref{2.1}, \eqref{2.2}, we introduce auxiliary vector variables, following the mixed finite element methodology \cite{boffi2013mixed}.
Define the flux vector as
\begin{equation}\label{2.4}
  \bm q(\bm x) = - K(\bm x)\, \operatorname{grad} u(\bm x).
\end{equation}
With this notation, equation \eqref{2.1} can be written as
\begin{equation}\label{2.5}
	\operatorname{div} \bm q + r(\bm x)\, u(\bm x) = f(\bm x), 
  \quad \bm x \in \Omega .
\end{equation}
In the finite-volume method, we approximate the system \eqref{2.4}--\eqref{2.5} together with the boundary condition \eqref{2.2}, using discrete analogues of vector operators.

The finite-volume construction involves the introduction and subsequent elimination of flux variables.  
The main stages are as follows:
\begin{itemize}
 \item Assign the primary variable $u(\bm x)$ on a mesh, seeking its approximation at the mesh nodes. These nodes define a partition of the computational domain into primary cells.
 \item Introduce a separate mesh for the auxiliary flux variable $\bm q(\bm x)$. Each primary mesh node is associated with a control volume, and fluxes are evaluated at the boundaries of these volumes.
 \item For each control volume, write a discrete analogue of equation \eqref{2.5} by integrating over the control volume and incorporating primary and auxiliary mesh quantities.
 \item Approximate the flux variables using differential relations such as \eqref{2.4}.
 \item Eliminate the auxiliary fluxes to obtain a discrete problem solely for the primary unknowns.
\end{itemize}

In this way, the finite volume method provides a constructive framework for developing mesh-based formulations of boundary value problems on irregular meshes.  
Below we consider the key elements of this approach for problem \eqref{2.1}, \eqref{2.2}.
\section{Meshes and mesh functions}

Meshes are constructed for the approximation of scalar and vector functions.
The primary mesh is a Delaunay triangulation without obtuse angles (D-mesh). 
Based on it, we build a Voronoi diagram (V-mesh). 
From both meshes we introduce a merged Voronoi--Delaunay mesh (MVD-mesh) with quadrilateral cells, in which vector function values are defined at the intersections of the diagonals.

\subsection{Delaunay mesh}

As the primary mesh, we consider a triangular mesh in the computational domain $\Omega$, obtained by Delaunay triangulation. 
We assume that $\Omega$ is a convex polygon with boundary $\partial \Omega$.

\begin{figure}[htbp]
\includegraphics[width=0.75\linewidth]{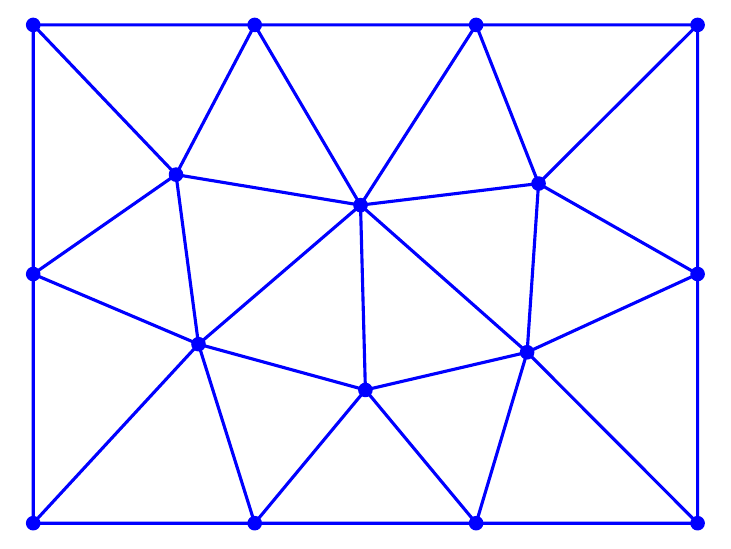}
\caption{Delaunay triangulation}
\label{f-1}
\end{figure}

The Delaunay triangulation is assumed to contain no obtuse angles. 
Such triangulations are often called non-obtuse meshes. 
A simple example is shown in Figure~\ref{f-1}.  
The usual mesh quality requirement \cite{cheng2013delaunay} is related to the minimum interior angles of the triangles. 
This characteristic is typically used in mesh generation algorithms and in software design. 
Computational algorithms for controlling the maximum angle of triangulations are described in \cite{bern2017triangulations}.

In the closure $\overline{\Omega} = \Omega \cup \partial \Omega$, we denote the D-mesh nodes as $\bm x_i^D$, $i=1,2,\ldots,M_D$:
\[
 \overline{\omega}^D = \{ \bm x \ | \ \bm x = \bm x_i^D, \  \bm x_i^D \in \overline{\Omega}, \ i = 1,2, \ldots, M_D \}.
\]
We distinguish interior nodes $\omega^D$ and boundary nodes $\partial \omega^D$, which include the vertices of $\Omega$:  
$\overline{\omega}^D = \omega^D \cup \partial \omega^D$.

\subsection{Voronoi mesh}

Each node $\bm x_i^D$, $i=1,2,\ldots,M_D$, is associated with a Voronoi polygon $\Omega_i^D$, or its intersection with $\Omega$. 
The Voronoi polygon for node $i$ is the set of points in $\Omega$ closer to $\bm x_i^D$ than to any other node $\bm x_{i^\pm}^D$:
\[
 \Omega_i^D = \{ \bm x \in \Omega \ | \ |\bm x - \bm x_i^D| < |\bm x - \bm x_{i^\pm}^D|, \ i^\pm=1,2,\ldots,M_D, \ i^\pm \neq i \}, 
 \quad i=1,2,\ldots,M_D .
\]

If the primary mesh is acute, then the Voronoi vertices inside $\Omega$ always lie within the corresponding triangles. 
Vertices of Voronoi polygons lying on $\partial \Omega$ are the midpoints of Delaunay edges belonging to the boundary.  

\begin{figure}[htbp]
\includegraphics[width=0.75\linewidth]{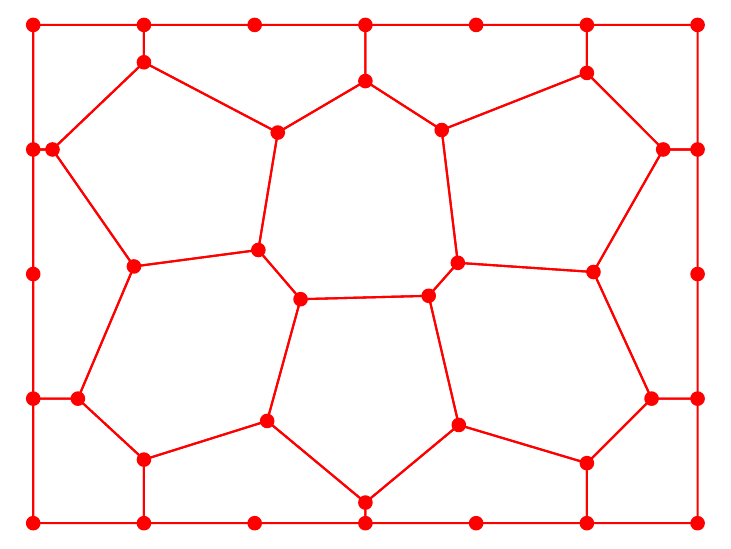}
\caption{Voronoi diagram}
\label{f-2}
\end{figure}

Let $\bm x_j^V$, $j=1,2,\ldots,M_V$, be the Voronoi vertices. 
For the triangulation in Figure~\ref{f-1}, the Voronoi mesh is shown in Figure~\ref{f-2}.

The set of all Voronoi vertices defines the V-mesh:
\[
 \overline{\omega}^V = \{ \bm x \ | \ \bm x = \bm x_j^V, \  \bm x_j^V \in \overline{\Omega}, \ j=1,2,\ldots,M_V \}, 
\]
with $\overline{\omega}^V = \omega^V \cup \partial \omega^V$.

In discretization, mesh cells are defined by their vertices, while in the finite volume method the domain is additionally partitioned into control volumes associated with mesh nodes.

The control volumes for the D-mesh are the regions $\Omega_i^D$:
\[
 \overline{\Omega} = \bigcup_{i=1}^{M_D} \overline{\Omega}_i^D, \quad 
 \overline{\Omega}_i^D = \Omega_i^D \cup \partial \Omega_i^D, \quad
 \Omega_i^D \cap \Omega_{i^\pm}^D = \emptyset, \ i^\pm \neq i .
\]

Similarly, each interior Voronoi node $\bm x_j^V$, $j=1,2,\ldots,M_V$, corresponds to a control volume $\Omega_j^V$, with
\[
 \overline{\Omega} = \bigcup_{j=1}^{M_V} \overline{\Omega}_j^V, \quad 
 \overline{\Omega}_j^V = \Omega_j^V \cup \partial \Omega_j^V, \quad
 \Omega_j^V \cap \Omega_{j^\pm}^V = \emptyset, \ j^\pm \neq j .
\]
For boundary nodes $\bm x_j^V \in \partial \Omega$, no control volumes are assigned: $\Omega_j^V = \emptyset$.

\subsection{Merged Voronoi--Delaunay mesh}

\begin{figure}[htbp]
 \includegraphics[width=0.75\linewidth]{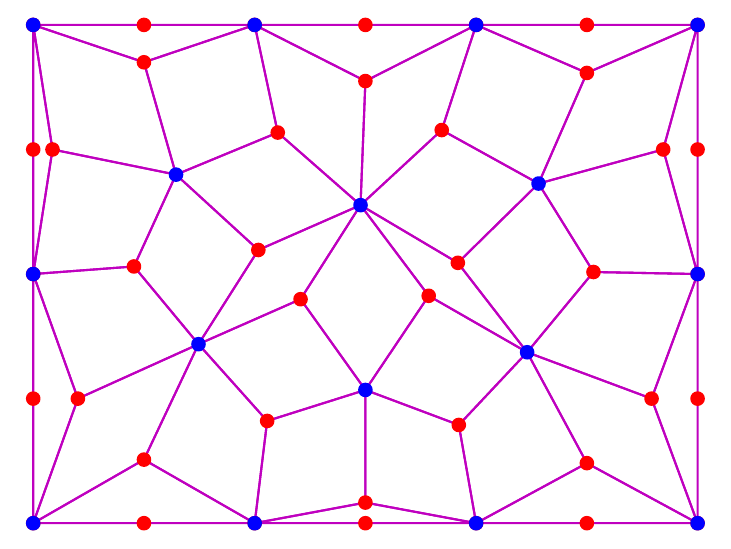}
 \caption{Merged Voronoi--Delaunay mesh}
 \label{f-3}
 \end{figure}

From the nodes of the D- and V-meshes we construct a partition of $\Omega$ into quadrilateral cells. 
Near the boundary these cells degenerate into triangles. 
The quadrilaterals are orthodiagonal, i.e., their diagonals are orthogonal.  
The center of each MVD-cell $\Omega_m$, $m=1,2,\ldots,M$, is the intersection of its diagonals.  
An example of the MVD-mesh, corresponding to Figures~\ref{f-1} and \ref{f-2}, is shown in Figure~\ref{f-3}.

\begin{figure}[htbp]
 \includegraphics[width=0.75\linewidth]{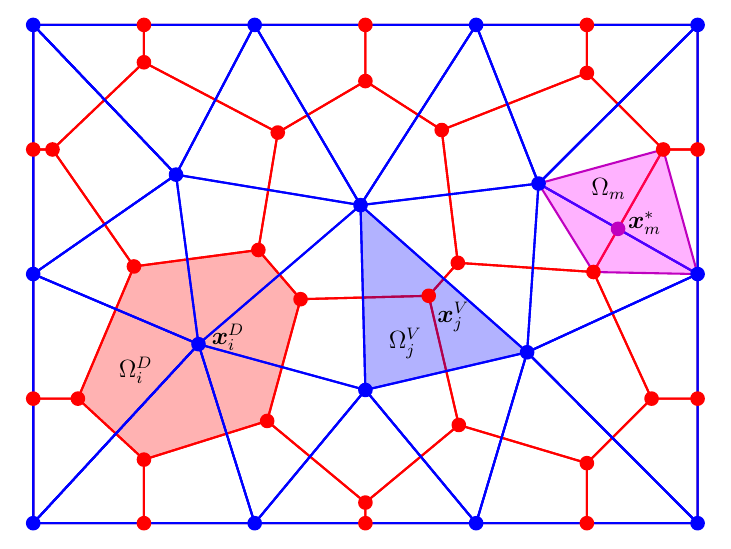}
 \caption{Cells of the D-, V-, and MVD-meshes}
 \label{f-4}
 \end{figure}

Figure~\ref{f-4} shows the three types of cells: Delaunay triangles, Voronoi polygons, and MVD orthodiagonal quadrilaterals.

\subsection{Scalar mesh functions}

For scalar functions $u(\bm x)$, $\bm x \in \Omega$, we introduced the Hilbert space $\mathcal{H} = L_2(\Omega)$ with its inner product and norm. 
Analogously, we define Hilbert spaces of mesh functions on the D- and V-meshes. 
Mesh functions are defined at mesh nodes, each associated with a control volume.

For D-mesh functions at nodes $\bm x_i^D$, $i=1,2,\ldots,M_D$, we introduce the space $H(\omega^D)$ with inner product
\[
 (y,v)_D = \sum_{\bm x \in \overline{\omega}^D} y(\bm x)\, v(\bm x)\, S^D(\bm x),
 \quad \|y\|_D = (y,y)_D^{1/2}.
\]
Here $S^D(\bm x)$ is the measure of the corresponding control volume:
\[
 S^D(\bm x) = \operatorname{meas}(\Omega_i^D), 
 \quad \bm x = \bm x_i^D \in \overline{\omega}^D .
\]
If homogeneous Dirichlet conditions \eqref{2.2} are imposed, then
\[
  H_0(\omega^D) = \{ u \in H(\omega^D) \ | \ u(\bm x)=0, \ \bm x \in \partial\omega^D \},
\]
and the inner product reduces to
\[
  (y,v)_D = \sum_{\bm x \in \omega^D} y(\bm x)\, v(\bm x)\, S^D(\bm x).
\]

Similarly, for V-mesh functions at interior nodes $\bm x_j^V$, $j=1,2,\ldots,M_V$, we define the space $H(\omega^V)$:
\[
 (y,v)_V = \sum_{\bm x \in \omega^V} y(\bm x)\, v(\bm x)\, S^V(\bm x),
 \quad \|y\|_V = (y,y)_V^{1/2},
\]
where
\[
 S^V(\bm x) = \operatorname{meas}(\Omega_j^V), 
 \quad \bm x = \bm x_j^V \in \omega^V .
\]

\subsection{Vector mesh functions}

The centers of MVD-cells are denoted $\bm x_m^*$, $m=1,2,\ldots,M$, and the mesh is $\omega^*$. 
Vector mesh functions are defined at these nodes.

\begin{figure}[htbp]
 \includegraphics[width=0.75\linewidth]{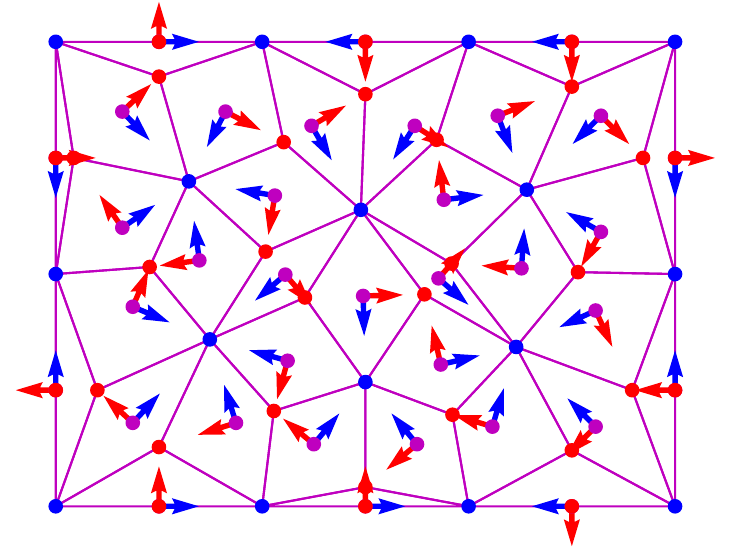}
 \caption{Local coordinate systems on the MVD-mesh}
 \label{f-5}
 \end{figure}

In each MVD-cell we introduce a local orthogonal coordinate system with unit vectors $\bm e_D(\bm x)$ and $\bm e_V(\bm x)$. 
The components of a vector $\bm v(\bm x) = \{v_D(\bm x), v_V(\bm x)\}$ are interpreted as tangential components on the faces of the Delaunay and Voronoi cells, respectively:
\begin{equation}\label{3.1}
  \bm v(\bm x) = \bm v_D(\bm x) + \bm v_V(\bm x),
  \quad \bm v_D(\bm x) = v_D(\bm x)\, \bm e_D(\bm x),
  \quad \bm v_V(\bm x) = v_V(\bm x)\, \bm e_V(\bm x),
  \quad \bm x_m^* \in \omega^* .
\end{equation}
The MVD-mesh with local coordinate directions is shown in Figure~\ref{f-5}.

For vector mesh functions at $\bm x_m^*$ we define the space $\bm H(\omega^*)$ with inner product
\[
 (\bm v,\bm w)_* = \sum_{\bm x \in \omega^*} \big( v_D(\bm x)\, w_D(\bm x) + v_V(\bm x)\, w_V(\bm x) \big)\, S^*(\bm x),
 \quad \|\bm v\|_* = (\bm v,\bm v)_*^{1/2}.
\]
Here
\[
 S^*(\bm x) = \operatorname{meas}(\Omega_m), 
 \quad \bm x = \bm x_m^* \in \omega^* .
\]

\section{Operator-difference approximations}

When studying boundary value problems for second-order elliptic and parabolic equations, the main focus is on approximations of the fundamental differential operators of vector calculus: the gradient of a scalar function and the divergence of a vector function.

\subsection{Mesh gradient operator}

A scalar mesh function $y(\bm x)$ is defined by its values at the nodes of the Delaunay mesh,
\[
  y_D(\bm x) = y(\bm x), \quad \bm x \in \overline{\omega}^D,
\]
and at the nodes of the Voronoi mesh,
\[
  y_V(\bm x) = y(\bm x), \quad \bm x \in \overline{\omega}^V .
\]

\begin{figure}[htbp]
\includegraphics[width=0.5\linewidth]{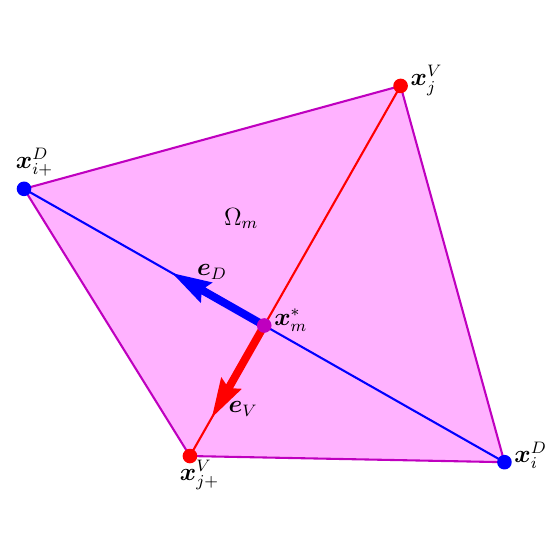}
\caption{Gradient approximation}
\label{f-6}
\end{figure}

In the local coordinate system (see Figure~\ref{f-6}), the mesh gradient is written as
\begin{equation}\label{4.1}
(\operatorname{grad}_h y)(\bm x) = (\operatorname{grad}_D y_D)(\bm x) + (\operatorname{grad}_V y_V)(\bm x),
\quad  \bm x = \bm x_{m}^* \in \omega^* .
\end{equation}
The components of the gradient vector are associated with finite differences along the faces of the D- and V-meshes.  
We set
\[
  (\operatorname{grad}_D y_D)(\bm x) = (\partial_D y_D) (\bm x) \, \bm e_D(\bm x),
  \quad
  (\operatorname{grad}_V y_V)(\bm x) = (\partial_V y_V) (\bm x) \, \bm e_V(\bm x),
  \quad \bm x \in \omega^* .
\]

Thus, for the components we obtain
\[
  \operatorname{grad}_D y_D(\bm x_{m}^*) = \frac{y(\bm x_{i^+}^D) - y(\bm x_{i}^D)}{|\bm x_{i^+}^D - \bm x_{i}^D|} \, \bm e_D(\bm x_{m}^*),
\]
\[
  \operatorname{grad}_V y_V(\bm x_{m}^*) = \frac{y(\bm x_{j^+}^V) - y(\bm x_{j}^V)}{|\bm x_{j^+}^V - \bm x_{j}^V|} \, \bm e_V(\bm x_{m}^*).
\]
Hence,
\[
  (\partial_D y_D) (\bm x_{m}^*) = \frac{y(\bm x_{i^+}^D) - y(\bm x_{i}^D)}{|\bm x_{i^+}^D - \bm x_{i}^D|},
  \quad
  (\partial_V y_V) (\bm x_{m}^*) = \frac{y(\bm x_{j^+}^V) - y(\bm x_{j}^V)}{|\bm x_{j^+}^V - \bm x_{j}^V|}.
\]

We distinguish subspaces of vector mesh functions $\bm H(\omega^*)$ by
\[
  \bm H(\omega^*) = \bm H_D(\omega^*) + \bm H_V(\omega^*),
\]
where
\[
  \bm H_D(\omega^*) = \{ \bm v \in \bm H(\omega^*) \mid \bm v_V(\bm x) = 0, \ \bm x \in \omega^*\},
\]
\[
  \bm H_V(\omega^*) = \{ \bm v \in \bm H(\omega^*) \mid \bm v_D(\bm x) = 0, \ \bm x \in \omega^*\}.
\]
Accordingly, the domains and ranges of the gradient operators are
\begin{equation}\label{4.2}
 \operatorname{grad}_D: \ H(\omega^D) \to \bm H_D(\omega^*),
 \quad
 \operatorname{grad}_V: \ H(\omega^V) \to \bm H_V(\omega^*).
\end{equation}

\subsection{Flux vector}

The flux vector is defined according to \eqref{2.4}.  
Within each MVD-cell, the diffusion tensor must be expressed in the corresponding local coordinate system.

Let $\theta(\bm x_m^*)$ denote the rotation angle from Cartesian to the local coordinates of cell $\Omega_m$.  
The clockwise rotation matrix is
\[
 Q(\bm x_m^*) =
 \begin{pmatrix}
 	\cos \theta(\bm x_m^*)  & \sin \theta(\bm x_m^*) \\
 	- \sin \theta(\bm x_m^*) & \cos \theta(\bm x_m^*)
 \end{pmatrix}.
\]
In local coordinates, the diffusion tensor is
\[
  \tilde{K} (\bm x_m^*) = Q (\bm x_m^*) \, K (\bm x_m^*) \, Q^*(\bm x_m^*) ,
\]
with components
\[
  \tilde{K} (\bm x_m^*) =
   \begin{pmatrix}
   	k_{DD}(\bm x_m^*) & k_{DV}(\bm x_m^*) \\
   	k_{VD}(\bm x_m^*) & k_{VV}(\bm x_m^*) 
   \end{pmatrix}.
\]

Taking into account \eqref{4.1}, the mesh flux vector is
\[
  \bm g (\bm x) = \bm g_D(\bm x) + \bm g_V(\bm x),
\quad  \bm x = \bm x_{m}^* \in \omega^* ,
\]
where
\begin{equation}\label{4.3}
\begin{split}
  g_D(\bm x) & = - k_{DD} (\bm x) (\partial_D y_D)(\bm x) - k_{DV}(\bm x) (\partial_V y_V)(\bm x), \\
  g_V(\bm x) & = - k_{VD} (\bm x) (\partial_D y_D)(\bm x) - k_{VV}(\bm x) (\partial_V y_V)(\bm x).
\end{split}
\end{equation}
Thus, in the local coordinate system we obtain convenient formulas for the approximate computation not only of the gradient but also of the flux vector in anisotropic media.

\subsection{Mesh divergence operator}

The divergence approximation is constructed separately at Delaunay and Voronoi mesh nodes when vector functions are given at $\omega^*$ (the centers of MVD-cells).  
The approach relies on the integral equality for the control volume.

For $\bm v(\bm x)$ defined at $\bm x \in \omega^*$, the approximation
\[
  y(\bm x) = \operatorname{div}_h \bm v(\bm x)
\]
is evaluated at $\bm x \in \omega$.  
At Delaunay nodes we write $\operatorname{div}_D$, and at Voronoi nodes — $\operatorname{div}_V$.

\begin{figure}[htbp] \includegraphics[width=0.5\linewidth]{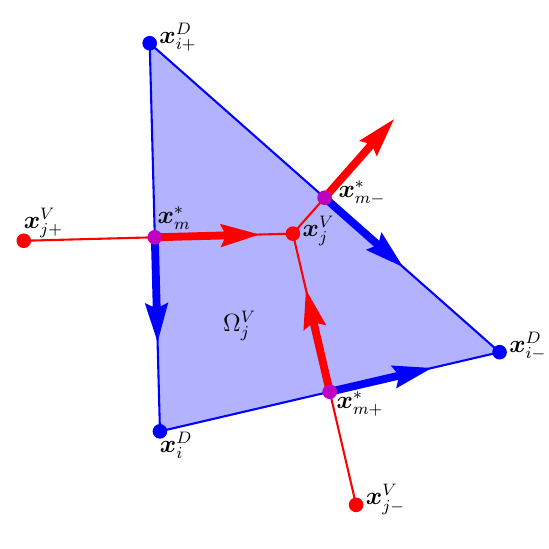} \caption{Divergence approximation at a V-mesh node} \label{f-7} \end{figure}

Consider the case where a node of $\omega$ belongs to the Voronoi diagram, so the control volume is a Delaunay cell (Figure~\ref{f-7}).  
By the divergence theorem,
\begin{equation}\label{4.4}
\int_{\Omega_j^V} \operatorname{div} \bm w (\bm x) \, d \bm x = \int_{\partial \Omega_j^V} \bm w \cdot \bm n \, d \bm x .
\end{equation}
Thus, the operator $\operatorname{div}_V: \bm H(\omega^*) \to H(\omega^V)$ is defined as
\begin{equation}\label{4.5}
(\operatorname{div}_V \bm v)(\bm x_j^V) = \frac{1}{S(\bm x_j^V)} \int_{\partial \Omega_j^V} \widetilde{\bm v} \cdot \bm n \, d \bm x ,
\end{equation}
where $\widetilde{\bm v}(\bm x)$ is the piecewise constant interpolation of $\bm v$ on $\partial \Omega_j^V$.

Let $\bm n_j^V(\bm x)$ be the outward normal to $\partial \Omega_j^V$.  
Using the notations of Figure~\ref{f-7}, from \eqref{4.5} we obtain
\begin{equation}\label{4.6}
\begin{split}
(\operatorname{div}_V \bm v)(\bm x_j^V) = \frac{1}{S(\bm x_j^V)} \Big(&
 v_V(\bm x_m^*) \delta_j^V(\bm x_m^*) |\bm x_{i^+}^D - \bm x_{i}^D| \\
&+ v_V(\bm x_{m^+}^*) \delta_j^V(\bm x_{m^+}^*) |\bm x_{i}^D - \bm x_{i^-}^D| \\
&+ v_V(\bm x_{m^-}^*) \delta_j^V (\bm x_{m^-}^*) |\bm x_{i^+}^D - \bm x_{i^-}^D| \Big).
\end{split}
\end{equation}
Here $\delta_j^V(\bm x) = \bm e_V (\bm x) \cdot \bm n_j^V(\bm x)$ equals $1$ if $\bm e_V (\bm x)$ coincides with $\bm n_j^V(\bm x)$, and $-1$ otherwise.

Similarly, at Delaunay nodes the operator $\operatorname{div}_D$ is constructed using expressions like
\[
\frac{1}{S(\bm x_i^D)} v_D(\bm x_{m}^*) \delta_i^D(\bm x_{m}^*) |\bm x_{j^+}^V - \bm x_{j}^V|,
\quad
\frac{1}{S(\bm x_i^D)} v_D(\bm x_{m^+}^*) \delta_i^D(\bm x_{m^+}^*) |\bm x_{j}^V - \bm x_{j^-}^V|,
\]
where $\delta_i^D(\bm x) = \bm e_D(\bm x) \cdot \bm n_i^D(\bm x)$, and $\bm n_i^D(\bm x)$ is the outward normal to the Voronoi control volume $\Omega_i^D$.

Finally, taking into account representation \eqref{3.1}, we write
\[
  (\operatorname{div}_h \bm v)(\bm x) =
  \begin{cases}
    (\operatorname{div}_D \bm v_D)(\bm x), & \bm x \in \omega^D, \\
    (\operatorname{div}_V \bm v_V)(\bm x), & \bm x \in \omega^V ,
  \end{cases}
\]
with
\begin{equation}\label{4.7}
  \operatorname{div}_D: \bm H_D(\omega^*) \to H(\omega^D),
  \quad
  \operatorname{div}_V: \bm H_V(\omega^*) \to H(\omega^V).
\end{equation}
Thus, $\operatorname{div}_D$ is defined on $\bm H_D(\omega^*)$ with values in $H(\omega^D)$, while $\operatorname{div}_V$ is defined on $\bm H_V(\omega^*)$ with values in $H(\omega^V)$.

\section{Mesh boundary value problem}

A discrete analogue of the boundary value problem \eqref{2.1}, \eqref{2.2} with sufficiently smooth coefficients and right-hand side is constructed using the finite volume method. To this end, the mixed formulation \eqref{2.2}, \eqref{2.4}, \eqref{2.5} is employed, together with the mesh gradient and divergence operators defined on the considered grids.

\subsection{Finite-volume approximation}

The discrete system is derived by approximating the gradient and divergence operators in equations \eqref{2.4}, \eqref{2.5}.  
Equation \eqref{2.5} leads to separate relations for the interior nodes of the Delaunay and Voronoi meshes:
\begin{equation}\label{5.1}
	\operatorname{div_D} \bm g_D + r_D(\bm x) y_D = f_D(\bm x), 
	\quad \bm x \in \omega^D ,
\end{equation}
\begin{equation}\label{5.2}
	\operatorname{div_V} \bm g_V + r_V(\bm x) y_V = f_V(\bm x), 
	\quad \bm x \in \omega^V .
\end{equation}

The approximation of \eqref{2.4} is performed according to \eqref{4.3}.  
Substituting the mesh flux vector \eqref{4.3} into \eqref{5.1}, \eqref{5.2} yields a system of equations for the approximate solution:
\begin{equation}\label{5.3}
	- \operatorname{div_D} \Big ( k_{DD} \operatorname{grad}_D y_D \Big )
    - \operatorname{div_D} \Big ( k_{DV} \operatorname{grad}_V y_V \Big ) 
	+ r_D(\bm x) y_D = f_D(\bm x), 
	\quad \bm x \in \omega^D ,
\end{equation}
\begin{equation}\label{5.4}
	- \operatorname{div_V} \Big ( k_{VD} \operatorname{grad}_D y_D \Big )
    - \operatorname{div_V} \Big ( k_{VV} \operatorname{grad}_V y_V \Big ) 
	+ r_V(\bm x) y_V = f_V(\bm x), 
	\quad \bm x \in \omega^V .
\end{equation}

For boundary nodes, the mesh functions satisfy homogeneous Dirichlet conditions:
\begin{equation}\label{5.5}
  y_D(\bm x) = 0, \quad \bm x \in \partial \omega^D, 
  \qquad 
  y_V(\bm x) = 0, \quad \bm x \in \partial \omega^V .
\end{equation}

The coupled nature of the system \eqref{5.3}, \eqref{5.4} appears whenever $k_{DV}(\bm x) = k_{VD}(\bm x) \neq 0$ for at least some nodes of the MVD mesh, $\bm x \in \omega^*$.  
In the particular case of a scalar diffusion coefficient, $K(\bm x) = k(\bm x)$, the system reduces to
\begin{equation}\label{5.6}
	- \operatorname{div_D} \big( k \operatorname{grad}_D y_D \big )
    + r_D(\bm x) y_D = f_D(\bm x), 
	\quad \bm x \in \omega^D ,
\end{equation}
\begin{equation}\label{5.7}
	- \operatorname{div_V} \big( k \operatorname{grad}_V y_V \big ) 
	+ r_V(\bm x) y_V = f_V(\bm x), 
	\quad \bm x \in \omega^V .
\end{equation}

In this case, the approximate solutions are computed independently on the Delaunay and Voronoi meshes.  
If no strategy is used to combine these two solutions (e.g., to improve the overall accuracy), the use of two meshes becomes of limited practical value.

\subsection{Consistency of gradient and divergence approximations}

For sufficiently smooth continuous functions, the following integral identity holds:
\[
  \int_{\Omega} u(\bm x) \operatorname{div} \bm v \, d \bm x +
  \int_{\Omega} \bm v(\bm x) \cdot \operatorname{grad} u \, d \bm x =
  \int_{\partial \Omega} u(\bm x) \bm v \cdot \bm n \, d \bm x .
\]
In the case $u(\bm x) = 0$ for $\bm x \in \partial \Omega$, this reduces to
\begin{equation}\label{5.8}
  \int_{\Omega} u(\bm x) \operatorname{div} \bm v \, d \bm x +
  \int_{\Omega} \bm v(\bm x) \cdot \operatorname{grad} u \, d \bm x = 0,
\end{equation}
which expresses that the operator adjoint to the divergence is equal to the gradient with a negative sign, i.e. $\operatorname{div}^* = - \operatorname{grad}$, in the corresponding Hilbert spaces.  
The discrete analogue of \eqref{5.8} is tied to the construction of consistent mesh approximations of the gradient and divergence, as realized in the framework of mimetic discretization methods \cite{shashkov2018conservative,da2014mimetic,vabishchevich2005finite}.

At the mesh level, the first term in \eqref{5.8} corresponds to scalar products in $H_0(\omega^D)$ and $H(\omega^V)$.  
To obtain a discrete analogue, we consider individual cells $\Omega_m$, $m=1,2,\ldots,M$, of the MVD mesh.  
For the Delaunay nodes, one has
\[
  I = \big( y_D, \operatorname{div}_D \bm v_D \big)_D, 
  \qquad 
  I = \sum_{m=1}^M I_m .
\]
When computing $I_m$ over the cell $\Omega_m$ (Figure~\ref{f-6}), two Delaunay nodes $\bm x_{i}^D$ and $\bm x_{i^+}^D$ are involved.  
In local coordinates, this gives
\[
\begin{split}
  I_m & = y(\bm x_{i}^D) v_D(\bm x_m^*) |\bm x_{j^+}^V - \bm x_{j}^V|
        - y(\bm x_{i^+}^D) v_D(\bm x_m^*) |\bm x_{j^+}^V - \bm x_{j}^V| \\
      & = - v_D(\bm x_m^*) |\bm x_{j^+}^V - \bm x_{j}^V| 
          \big(y(\bm x_{i^+}^D) - y(\bm x_{i}^D)\big) \\
      & = - 2 v_D(\bm x_m^*) 
          \frac{y(\bm x_{i^+}^D) - y(\bm x_{i}^D)}{|\bm x_{i^+}^D - \bm x_{i}^D|} 
          S^*(\bm x_m^*) .
\end{split}
\]

For functions satisfying $y_D(\bm x)=0$ on $\partial \omega^D$, this yields a mesh analogue of \eqref{5.8}:
\[
  \big( y_D, \operatorname{div}_D \bm v_D \big)_D +
  2 \big(\bm v_D, \operatorname{grad}_D y_D \big)_* = 0.
\]
A similar argument applies for the Voronoi mesh, leading to
\[
  \big( y_V, \operatorname{div}_V \bm v_V \big)_V +
  2 \big(\bm v_V, \operatorname{grad}_V y_V \big)_* = 0.
\]
Thus, for the introduced function spaces,
\begin{equation}\label{5.9}
  \operatorname{div}^*_D = - 2 \operatorname{grad}_D,
  \qquad 
  \operatorname{div}^*_V = - 2 \operatorname{grad}_V .
\end{equation}

\subsection{A priori estimate}

For the mesh problem \eqref{5.3}--\eqref{5.5}, the standard issues of approximation, stability, and convergence of the discrete solution to the exact one are considered.  
Here we restrict ourselves to establishing a mesh analogue of the a priori estimate \eqref{2.3} for the boundary value problem \eqref{2.1}, \eqref{2.2}.

Multiply \eqref{5.3} scalarly in $H(\omega^D)$ by $y^D$, and \eqref{5.4} scalarly in $H(\omega^V)$ by $y^V$.  
Using \eqref{5.9} and the introduced notation, we obtain
\[
 2 \big( k_{DD}\operatorname{\partial}_D y_D, \operatorname{\partial}_D y_D \big)_* +
 2 \big( k_{DV}\operatorname{\partial}_V y_V, \operatorname{\partial}_D y_D \big)_* 
 + (r_D y_D, y_D)_D = (f_D, y_D)_D ,
\]
\[
 2 \big( k_{VD}\operatorname{\partial}_D y_D, \operatorname{\partial}_V y_V \big)_* +
 2 \big( k_{VV}\operatorname{\partial}_V y_V, \operatorname{\partial}_V y_V \big)_* 
 + (r_V y_V, y_V)_V = (f_V, y_V)_V .
\]

Adding these relations and using the properties of the coefficients in \eqref{2.1}, we arrive at
\begin{equation}\label{5.10}
  2 \kappa \Big( (\operatorname{\partial}_D y_D, \operatorname{\partial}_D y_D)_*
  + (\operatorname{\partial}_V y_V, \operatorname{\partial}_V y_V)_* \Big) 
  \leq (f_D, y_D)_D + (f_V, y_V)_V .
\end{equation}

For mesh functions on Delaunay and Voronoi grids, discrete Friedrichs-type inequalities hold.  
In \cite{SamVab2000gvm}, it is shown that
\[
  \|y_D\|_D \leq c_h \|\operatorname{\partial}_D y_D\|_* ,
\]
where the constant $c_h$ depends only on the computational domain $\Omega$.  
An analogous inequality holds for functions on the Voronoi mesh:
\[
  \|y_V\|_V \leq c_h \|\operatorname{\partial}_V y_V\|_* .
\]

Introducing the combined norm
\[
 \| \operatorname{grad}_h y\|_*^2 
 = \|\operatorname{\partial}_D y_D\|_*^2 + \|\operatorname{\partial}_V y_V\|_*^2 ,
\]
from \eqref{5.10} we obtain
\begin{equation}\label{5.11}
  \| \operatorname{grad}_h y\|_* 
  \leq \frac{c_h}{\kappa} \max \{ \|f_D\|_D, \|f_V\|_V \} .
\end{equation}

Inequality \eqref{5.11} is the mesh analogue of the continuous a priori estimate \eqref{2.3}.
\section{Numerical experiments}

We present the results of the approximate solution of a diffusion–reaction problem in a rectangle using the finite volume method on a merged Voronoi–Delaunay mesh.
The main focus is on the accuracy of the approximate solution in various norms for a problem with a known solution under changing sizes of the computational meshes.

\subsection{Test problem}

Below are the results of calculations for the boundary value problem \eqref{2.1}, \eqref{2.2} in the rectangle
\[
	\Omega = \{\bm x \  | \  \bm x = (x_1, x_2), \ 0 \le x_1 \le 1, \ 0 \le x_2 \le 0.75\} .
\]
We limited ourselves to the case of a homogeneous medium, where the diffusion tensor $K$ is constant. The main attention was on problems in anisotropic media, for which the developed computational technology of joint use of Delaunay and Voronoi meshes is most interesting.

In the calculations, the constant diffusion tensor $K$ was set to
\[
	K_1 =
	\begin{pmatrix}
		1 & 0 \\
		0 & 1 \\
	\end{pmatrix}, \quad
	K_2 =
	\begin{pmatrix}
		1 & 0 \\
		0 & 100 \\
	\end{pmatrix}, \quad
	K_3 =
	\begin{pmatrix}
		1 & 9 \\
		9 & 100 \\
	\end{pmatrix}, \quad
	K_4 =
	\begin{pmatrix}
		1 & -9 \\
		-9 & 100 \\
	\end{pmatrix}.
\]
Thus, problems with strong anisotropy of the medium were considered.
The right-hand side of \eqref{2.1} corresponded to the exact solution
\[
	u(\bm x) = x_1(1-x_1) \sin{\left (\frac43 \pi x_2 \right )}.
\]

The accuracy of the approximate solution $y(\bm x)$ was evaluated separately on the Delaunay mesh and on the Voronoi mesh. In the $L_2$ norm, we have
\[
\begin{split}
\varepsilon_2^D & = \|y - u\|_D,
\quad \|y - u\|_D^2 = \sum_{\bm x \in \omega^D} (y(\bm x) - u(\bm x))^2 S^D(\bm x),	\\
 \varepsilon_2^V & = \|y - u\|_V,
\quad \|y - u\|_V^2 = \sum_{\bm x \in \omega^V} (y(\bm x) - u(\bm x))^2 S^V(\bm x) ,
\end{split}
\]
and in the $L_\infty$ norm,
\[
\begin{split}
\varepsilon_\infty^D & = \max_{1 \leq i \leq M_D} |y(\bm x_i^D) - u(\bm x_i^D)|, \\
\varepsilon_\infty^V & = \max_{1 \leq j \leq M_V} |y(\bm x_j^V) - u(\bm x_j^V)|.
\end{split}
\]
The accuracy of gradient computation at the nodes of the MVD-mesh was also evaluated:
\[
  \varepsilon = \|\operatorname{grad}_h y - \operatorname{grad} u\|_* .
\]

\subsection{Meshes}

A sequence of successively refined meshes was used for numerical experiments.
Triangular meshes were generated using the gmsh mesh generator \cite{geuzaine2009gmsh}, and based on them, Voronoi and MVD meshes were constructed.
No special efforts were made to generate meshes with guaranteed non-obtuse angles.

\begin{table}[h]
\caption{Mesh characteristics.}
\begin{center}
\begin{tabular}{|c|r|c|r|r|}
\hline
Mesh & D-mesh ($M_D$) & min/max angle & V-mesh ($M_V$) & MVD-mesh ($M$) \\ 
\hline
1 & 16 & 42.7 / 81.2 & 30 & 35 \\ 
4 & 36 & 42.9 / 83.3 & 70 & 87 \\ 
7 & 129 & 42.2 / 85.3 & 256 & 346 \\ 
10 & 427 & 43.9 / 87.2 & 852 & 1204 \\ 
13 & 1266 & 42.9 / 85.5 & 2530 & 3665 \\ 
16 & 4432 & 40.3 / 85.4 & 8862 & 13047 \\ 
19 & 16879 & 40.8 / 85.4 & 33756 & 50150 \\ 
22 & 64601 & 41.3 / 85.3 & 129200 & 192848 \\ 
\hline
\end{tabular}
\end{center}
\label{table}
\end{table}

We used 22 meshes in total.
The starting (coarsest) meshes are shown in Figs.~\ref{f-1}--\ref{f-3}.
The number of nodes for three meshes is given in Table~\ref{table}.
The quality of the initial Delaunay triangulation can be assessed from the data on the minimum and maximum angles of the triangular mesh cells.

\subsection{Accuracy of the approximate solution}

\begin{figure}[h]
\includegraphics[width=0.75\textwidth]{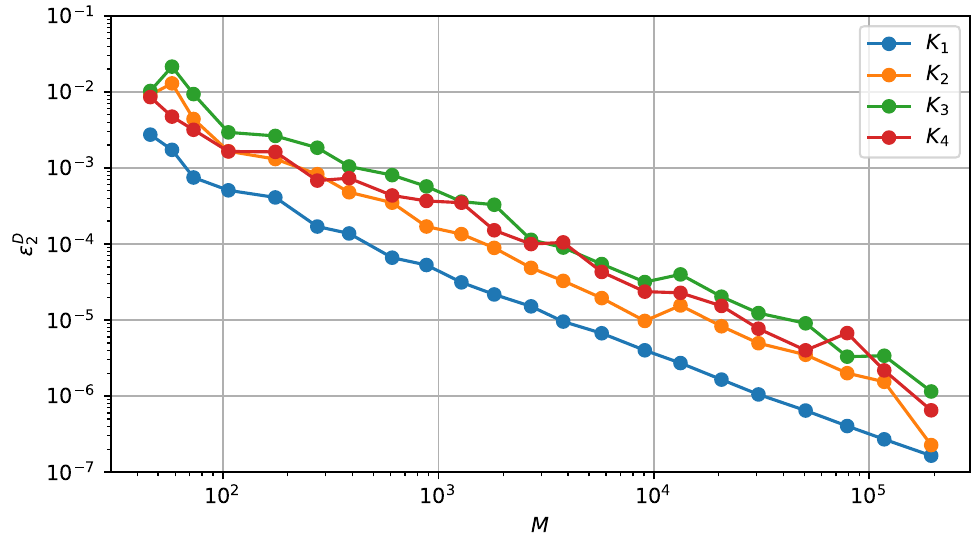}
\caption{Error of the approximate solution on the Delaunay mesh in the $L_2$ norm.}
\label{f-8}
\end{figure}

\begin{figure}[h]
\includegraphics[width=0.75\textwidth]{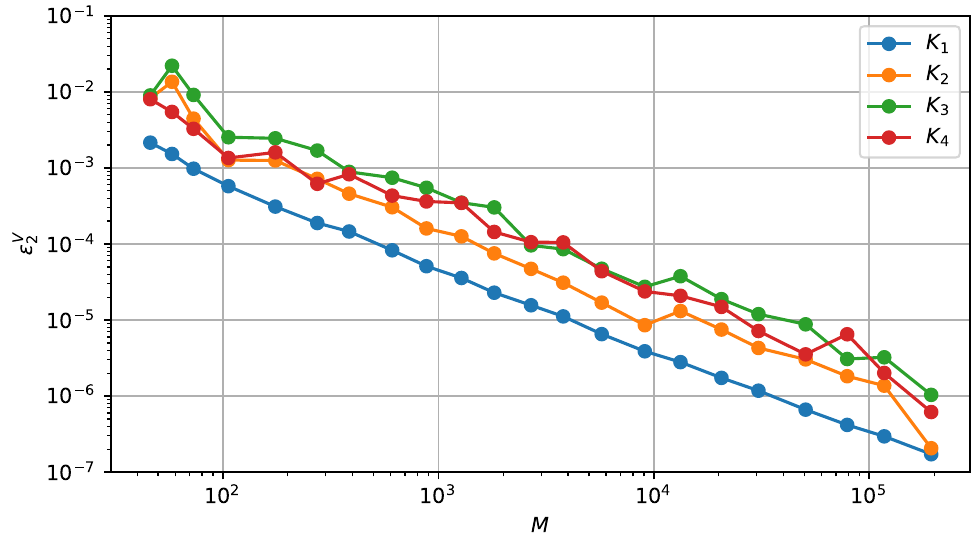}
\caption{Error of the approximate solution on the Voronoi mesh in the $L_2$ norm.}
\label{f-9}
\end{figure}

\begin{figure}[h]
\includegraphics[width=0.75\textwidth]{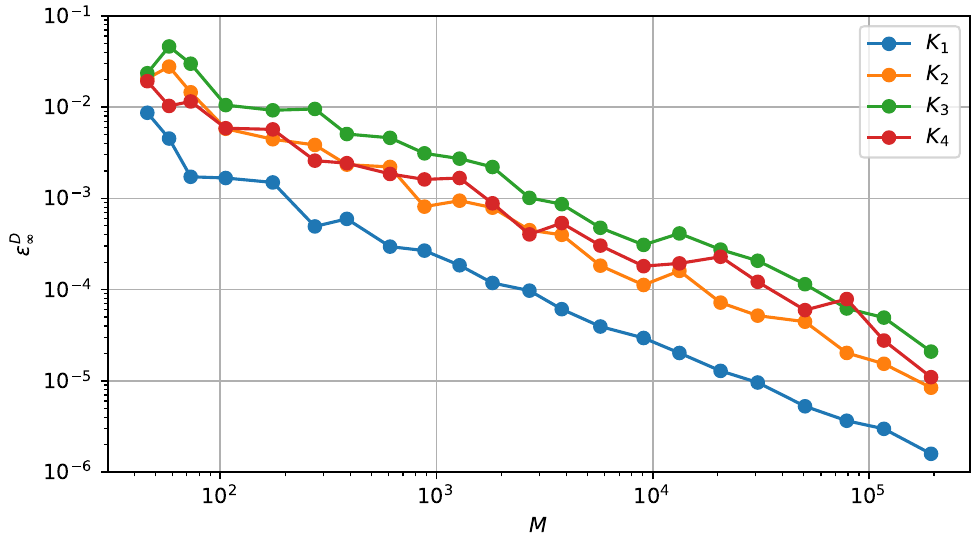}
\caption{Error of the approximate solution on the Delaunay mesh in the $L_\infty$ norm.}
\label{f-10}
\end{figure}

\begin{figure}[h]
\includegraphics[width=0.75\textwidth]{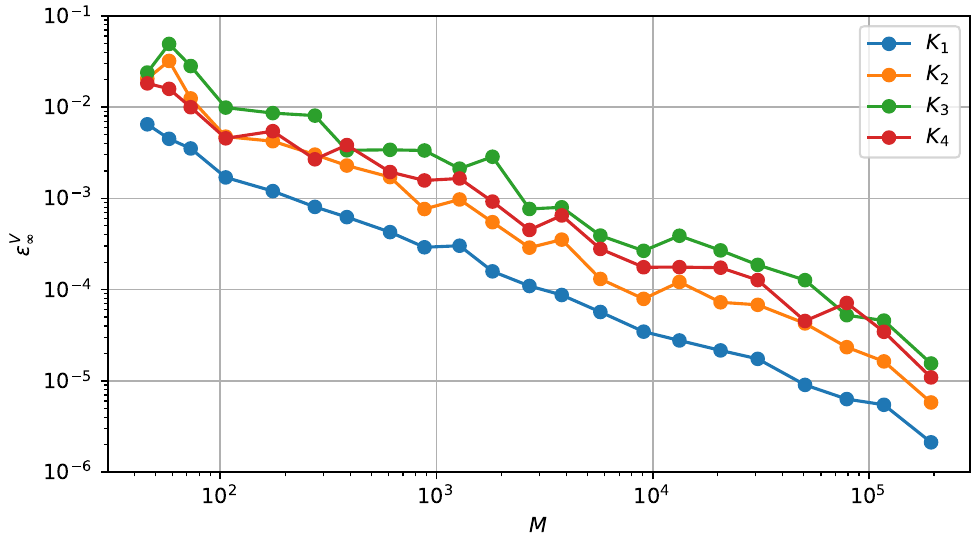}
\caption{Error of the approximate solution on the Voronoi mesh in the $L_\infty$ norm.}
\label{f-11}
\end{figure}

\begin{figure}[h]
\includegraphics[width=0.75\textwidth]{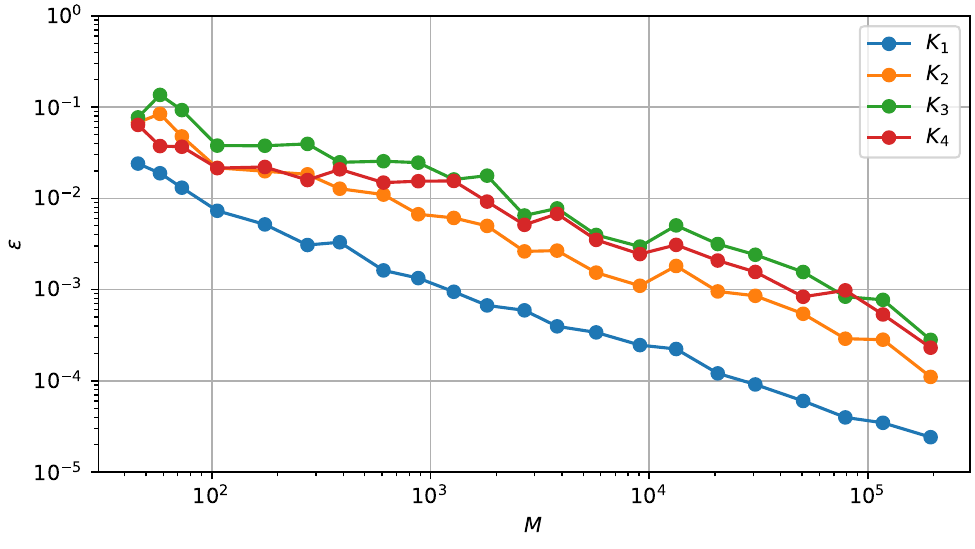}
\caption{Error of gradient computation.}
\label{f-12}
\end{figure}

The error norms of the approximate solution for different diffusion tensors are shown in Figs.~\ref{f-8}--\ref{f-12}.
The main conclusions can be formulated as follows:

\begin{itemize}

\item
The error of the approximate solution in the $L_2$ norm decreases approximately linearly with the number of nodes on both the Delaunay mesh and the Voronoi mesh (Figs.~\ref{f-8}, \ref{f-9}).

\item
In the uniform $L_\infty$ norm, the error decreases, preserving the same asymptotic dependence on the number of mesh nodes (Figs.~\ref{f-10}, \ref{f-11}).

\item
The error of gradient computation remains sufficiently small (Fig.~\ref{f-12}).

\item
The dependence on medium anisotropy is pronounced.
We observe a significant drop in accuracy when solving the diffusion problem in an anisotropic medium compared to the isotropic case.
This is most evident in the approximate computation of the gradient.

\end{itemize}

\section{Conclusions}

\begin{enumerate}[(1)]
    \item
    For a Delaunay triangulation with acute-angled cells, a Voronoi diagram is constructed.  
    In computational practice, the mesh solution is usually determined at Delaunay nodes in two variants of domain partitioning:  
    (a) node-centered approximation using Delaunay triangles, and  
    (b) cell-centered approximation using Voronoi polygons.  
    Merged Voronoi--Delaunay (MVD) meshes are introduced, whose nodes include both Delaunay vertices and Voronoi vertices.  
    The cells of such a mesh are orthodiagonal quadrilaterals.

    \item
    In the finite volume method, both the solution and fluxes are approximated.  
    The approximate solution is obtained by integrating the equation over a control volume.  
    In the MVD approach, for Delaunay nodes the control volumes are Voronoi cells, while for Voronoi nodes the control volumes are Delaunay cells.  
    Fluxes are approximated at the centers of MVD cells, i.e.\ at the intersections of the orthogonal diagonals.

    \item
    On irregular computational meshes, accurate discretization requires approximations of invariant vector (tensor) analysis operators.  
    Discrete analogues of gradient and divergence operators have been constructed on MVD meshes in a form consistent with the Mimetic Discretization methodology.

    \item
    The developed computational technology has been applied to a two-dimensional diffusion--reaction problem in anisotropic media.  
    The diffusion coefficient is represented by a symmetric tensor of rank two.  
    A coupled system of equations is formulated for the mesh solution at both Delaunay and Voronoi nodes.  
    Numerical experiments confirm the accuracy, robustness, and efficiency of the proposed approach.
\end{enumerate}


\end{document}